\documentclass[10pt,leqno]{article}
\usepackage{latexsym}
\usepackage{amsthm}
\usepackage{amsfonts}
\usepackage{amssymb}
\usepackage{amsmath} 
\usepackage{graphicx}
\usepackage{url} 
\usepackage{color}
\usepackage{indentfirst}

\usepackage[top=30truemm,bottom=30truemm,left=30truemm,right=30truemm]{geometry}

\theoremstyle{plain}
  \newtheorem{theorem}{Theorem}[section] 
  
  \newtheorem{proposition}{Proposition}[section]
  \newtheorem{lemma}{Lemma}[section]

\theoremstyle{remark}
  \newtheorem{remark}{Remark}[section]

\theoremstyle{definition}
  
  \newtheorem{notation}{Notation}[section]

\include{psfig}

\makeatletter

\@addtoreset{equation}{section}
\makeatother

\include{psfig}

\begin{document}

\markboth{Hideo Takaoka}{Derivative nonlinear Schr\"odinger equations}

\title{Remarks on blow-up criteria for the derivative nonlinear Schr\"{o}dinger equation under the optimal threshold setting}
\author{Hideo Takaoka\thanks{This work was supported by JSPS KAKENHI Grant Number 18H01129.}\\
Department of Mathematics, Kobe University\\
Kobe, 657-8501, Japan\\
takaoka@math.kobe-u.ac.jp}

\date{\empty}

\maketitle

\begin{abstract}
We study the Cauchy problem of the mass critical nonlinear Schr\"odinger equation with derivative with the $4\pi$ mass.
One has the global well-posedness in $H^1$ whenever ``the mass is strictly less than $4\pi$" or whenever ``the mass is equal to $4\pi$ and the momentum is strictly less than zero".
In this paper, by the concentration compact principle as originally done by Kenig-Merle, we obtain the limiting profile of blow up solutions with the critical $4\pi$ mass.
\end{abstract}

{\it $2010$ Mathematics Subject Classification Numbers.}
35Q55, 42B37.


\section{Introduction}\label{sec:introduction}

In this paper, we consider the Cauchy problem of the derivative nonlinear Schr\"odinger equation;
\begin{eqnarray}\label{dnls}
\left\{
\begin{array}{ll}
i\partial_t u+\partial_{x}^2u=i\partial_x(|u|^2u), & (t,x)\in\mathbb{R}^2,\\
u(0,x)=u_0(x), & x\in \mathbb{R},
\end{array}
\right.
\end{eqnarray}
where $u$ is a complex valued function of $(t,x)\in \mathbb{R}^2$.
Our aim is to provide the profile of solutions to the Cauchy problem (\ref{dnls}) with $\|u_0\|_{L^2}^2=4\pi$, if the solution of the current equation may became unbounded with respect to the energy norm $H^1(\mathbb{R})$.

The Cauchy problem (\ref{dnls}) is $L^2$ scale critical in the sense that the $L^2$ norm is invariant under the scaling transform;
\begin{eqnarray*}
u_{\lambda}(t,x)=\lambda^{1/2}u(\lambda^2 t,\lambda x),\quad \lambda>0.
\end{eqnarray*}
The equation in (\ref{dnls}) has an infinite family of conservations laws and can be solved by the inverse scattering method that are used to solve the Cauchy problem for certain classes of initial data \cite{kn}.
Formally solutions of (\ref{dnls}) satisfy the conservation laws, that, if $u=u(t,x)$ is solution to (\ref{dnls}), then for all $t\in\mathbb{R}$, the charge ($L^2$-norm)
\begin{eqnarray}\label{eq:mass}
M(u)=\|u\|_{L^2}^2=M(u_0),
\end{eqnarray}
the energy
\begin{eqnarray}\label{eq:H1energy}
E_1(u)=\|\partial_x u\|_{L^2}^2-\frac32\Im\langle \partial_x u,|u|^2u\rangle+\frac12\|u\|_{L^6}^6=E_1(u_0)
\end{eqnarray}
and the momentum
\begin{eqnarray}\label{eq:H1momentum}
P_1(u)=\Im\langle \partial_x u,u\rangle-\frac12\|u\|_{L^4}^4,
\end{eqnarray}
where $\langle \cdot,\cdot\rangle$ is the standard inner product in $L^2(\mathbb{R})$.
For the formal derivation of the recurrence formula of the infinite family of other conserved quantities, we refer the reader to the references \cite{kn,tsm1} (e.g. the conserved quantities $C_n$ in \cite{kn} and $Z^{(n)}$ in \cite{tsm1}).
As we will see below, next three specific conservation laws are formulated as the functional in the $H^3$ space.
Notice that two of those are the following potential Hamiltonians
\begin{eqnarray}\label{eq:H2energy}
E_2(u) & = & \|\partial_x^2u\|_{L^2}^2+\frac78\|u\|_{L^{10}}^{10}+\frac{25}{2}\left\||u|^2\partial_x u\right\|_{L^2}^2\\
& & +5\Re\langle (\partial_x u)^2,|u|^2u^2\rangle-5\Im\langle \partial_x^2u,|u|^2\partial_x u\rangle-\frac{35}{8}\Im\langle \partial_x u,|u|^6u\rangle,\nonumber\\
& = & E_2(u_0) \nonumber
\end{eqnarray}
and
\begin{eqnarray}\label{eq:H2momentum}
P_2(u) & = & \frac12\Im\langle \partial_x^2 u,\partial_x u\rangle-\frac{5}{16}\|u\|_{L^{8}}^8-2\|u\partial_xu\|_{L^2}^2\\
& & -\frac12\Re\langle (\partial_x u)^2,u^2\rangle+\frac54\Im\langle \partial_x u,|u|^4u\rangle\nonumber\\
& = & P_2(u_0).\nonumber
\end{eqnarray}
The next conservation laws involving third order derivative are
\begin{eqnarray}\label{eq:H3energy}
E_3(u) & = & \|\partial_x^3u\|_{L^2}^2+e_3(\partial_x^3u,\partial_x^2u,\partial_xu,u)\\
& = & E_3(u_0)\nonumber
\end{eqnarray}
and
\begin{eqnarray}\label{eq:H3momentum}
P_3(u) & = & \frac12\Im\langle \partial_x^3u,\partial_x^2u\rangle+p_3(\partial_x^2u,\partial_xu,u)\\
& = & P_3(u_0),\nonumber
\end{eqnarray}
where the functional $e_3(\partial_x^3u,\partial_x^2u,\partial_xu,u)$ and $p_3(\partial_x^2u,\partial_xu,u)$ satisfy
$$
|e_3(\partial_x^3u,\partial_x^2u,\partial_xu,u)|\le \varepsilon\|\partial_x^3u\|_{L^2}^2+C(\|u\|_{H^2},\varepsilon)
$$
for any $\varepsilon>0$, and
$$
|p_3(\partial_x^2u,\partial_xu,u)|\le C(\|u\|_{H^2}).
$$
Similarly, we have the conservation law at the forth order regularity
\begin{eqnarray}\label{eq:H4energy}
E_4(u) & = & \|\partial_x^4u\|_{L^2}^2+e_4(\partial_x^4u,\partial_x^3u,\partial_x^2u,\partial_xu,u)\\
& = & E_4(u_0),\nonumber
\end{eqnarray}
where the functional $e_4(\partial_x^4u,\partial_x^3u,\partial_x^2u,\partial_xu,u)$ satisfies
$$
|e_4(\partial_x^4u,\partial_x^3u,\partial_x^2u,\partial_xu,u)|\le \varepsilon\|\partial_x^4u\|_{L^2}^2+C(\|u\|_{H^3},\varepsilon)
$$
for any $\varepsilon>0$.

Let $\phi_{\omega,c}$ be a ground state solution to the following equation
\begin{eqnarray}\label{eq:elliptic}
-\phi''+\left(\omega-\frac{c^2}{4}\right)\phi+\frac{c}{2}\phi^3-\frac{3}{16}\phi^5=0.
\end{eqnarray}
The equation of (\ref{dnls}) possesses a two-parameter family of standing wave solution of the form
$$
q_{\omega,c}(t,x)=\phi_{\omega,c}(x+ct)e^{i\omega t-i\frac{c}{2}(x+ct)+i\frac34\int_{-\infty}^{x+ct}|\phi_{\omega,c}(y)|^2\,dy},
$$
where $c>0,~\omega\ge c^2/4$.
In the case when $\omega>c^2/4$, $\phi_{\omega,c}$ is an exponential decay function as $|x|\to\infty$ and
$$
\|\phi_{\omega,c}\|_{L^2}=\left(8\tan^{-1}\sqrt{\frac{\sqrt{4\omega}+c}{\sqrt{4\omega}-c}}\right)^{1/2}<\sqrt{4\pi}.
$$
For the borderline case $\omega=c^2/4$, the ground state of (\ref{eq:elliptic}) is expressed as
\begin{eqnarray}\label{eq:soliton}
Q_c(x)=\phi_{\frac{c^2}{4},c}(x)=\sqrt{\frac{4c}{(cx)^2+1}},
\end{eqnarray}
which satisfies $\|Q_c\|_{L^2}^2=4\pi$ for any $c>0$. 
It is worth noting that the standing wave solutions satisfy $E_1(q_{\omega,c})=-c\sqrt{4\omega-c^2}$ and $P_1(q_{\omega,c})=2\sqrt{4\omega-c^2}$, that imply $E_1(q_{c^2/4,c})=P_1(q_{c^2/4,c})=0$.
Moreover, $E_2(q_{c^2/4,c})=P_2(q_{c^2/4,c})=P_3(q_{c^2/4,c})=0$ (see Section \ref{sec:P2Q}).
We infer that $E_k(q_{c^2/4,c})=P_k(q_{c^2/4,c})=0$ for any $k\in\mathbb{N}$, but the author does not know whether it is true.

The well-posedness of the Cauchy problem (\ref{dnls}) has been extensively studied.
In \cite{h1,h2}, Hayashi-Ozawa proved the local well-posedness in the energy space $H^1(\mathbb{R})$ and the global well-posedness if the initial data satisfy $\|u_0\|_{L^2}^2<2\pi$.
For rough data, it was proved that the local well-posedness holds in $H^s(\mathbb{R})$ for $s\ge 1/2$ \cite{t1} (see also \cite{csktt1} for the global solutions for rough initial data).
Recently, Wu \cite{w1} and Guo-Wu \cite{gw} obtained the global well-posedness above the mass threshold $2\pi$, more precisely, if the mass is less than $4\pi$, namely $\|u_0\|_{L^2}^2<4\pi$.
The main ingredient in the proof is the following result.

\begin{proposition}[Lemma 2.1 in \cite{gw}]\label{prop:gw}
If $\phi\in H^1(\mathbb{R})$ and $\phi\ne 0$, then
\begin{eqnarray*}
\int_{\mathbb{R}}\left(\Im(\overline{\phi}\partial_x\phi)+\frac14|\phi|^4\right)\,dx
\ge \frac14\|\phi\|_{L^4}^4\left(1-\frac{\|\phi\|_{L^2}}{\sqrt{4\pi}}\right)-\frac{4\sqrt{\pi}\|\phi\|_{L^2}}{\|\phi\|_{L^4}^4}\int_{\mathbb{R}}\left(|\partial_x\phi|^2-\frac{1}{16}|\phi|^6\right)\,dx.
\end{eqnarray*}
\end{proposition}
The estimate in Proposition \ref{prop:gw} has impact on the availability of the a priori estimate of solutions. 
In particular, using Proposition \ref{prop:gw}, Fukaya-Hayashi-Inui \cite{fhi} proved that the global well-posedness was extended to the initial data $\|u_0\|_{L^2}^2=4\pi$ and $P_1(u_0)<0$.
If $\|u_0\|_{L^2}^2=4\pi$ and $E_1(u_0)=P_1(u_0)=0$, then the solution of (\ref{dnls}) is the solitary wave solution and exists global in time (see \cite{fhi}).
More recently, Jenkins-Liu-Perry-Sulem \cite{jlps} proved the global well-posedness for the initial data in $H^{2,2}(\mathbb{R})$, by using the inverse scattering tools, where $H^{s,a}(\mathbb{R})$ is a weighted Sobolev space.
Note that soliton solutions $q_{c^2/4,c}$ do not belong to $H^{2,2}$.
Under the mass condition $\|u_0\|_{L^2}^2=4\pi$, it turned out that the global well-posedness for the initial data satisfying $P_1(u_0)\ge 0$ is open.

It is interesting to take the initial data $u_0\in H^1$ with $\|u_0\|_{L^2}=\sqrt{4\pi}$ and $P_1(u_0)\ge 0$.
In the following example, we will illustrate that the set of objects $u_0\in H^1$ that satisfy
\begin{eqnarray}\label{eq:example}
\|u_0\|_{L^2}=\sqrt{4\pi},\quad P_1(u_0)>0
\end{eqnarray}
is not empty.

\begin{proposition}\label{prop:exam}
There exists an element $u_0\in H^1$ satisfying (\ref{eq:example}).
\end{proposition}

This aim of this paper is to prove the blow-up phenomena with the asymptotic profile for the mass critical case.
The main result of this paper is the following theorem.

\begin{theorem}\label{thm:main}
Let $u_0\in H^1$ satisfy $\|u_0\|_{L^2}^2=4\pi$.
Assume that $u(t)$ be the corresponding solution of (\ref{dnls}) blows up in finite or infinite time, i.e. there exists $0<T\le +\infty$ such that
\begin{eqnarray}\label{eq:assum}
\limsup_{t\uparrow T}\|\partial_xu(t)\|_{L^2}=+\infty.
\end{eqnarray}
Then there exist a sequence denoted $\{R_n\}_{n\in\mathbb{N}}\subset H^1$ with $\gamma\in \mathbb{R},~t_n\in\mathbb{R},~x_n\in\mathbb{R}$ such that
$$
\lim_{n\to\infty}t_n=T
$$
and
\begin{eqnarray}\label{eq:bu}
u(t_n,x)=v_n(x)\exp\left(i\frac{3}{4}\int_{-\infty}^x|v_n(y)|^2\,dy \right),
\end{eqnarray}
where
$$
\lambda_n^{1/2}v_n(\lambda_n x)e^{\frac{i}{2}x}=e^{i\gamma}\left(Q_1(x-x_n)+R_n(x)\right),
$$
with 
$$
\lambda_n=\sqrt{\frac{3\cdot 2^3\cdot \pi}{\|v_n\|_{L^6}^6}},
$$
$$
\lim_{n\to\infty}\|\partial_xR_n\|_{L^2}=0
$$
and
$$
\lim_{n\to\infty}\|R_n\|_{L^p}=0
$$
for every $p\in(2,\infty)$.
\end{theorem}

\begin{remark}
Our assumption in Theorem \ref{thm:main} allows us to deduce
\begin{eqnarray}\label{eq:assum}
\lim_{n\to\infty}\lambda_n=0.
\end{eqnarray}
Indeed, let $u_0\in H^1$ with $P_1(u_0)\ge 0$ and assume the corresponding solution $u(t)$ of (\ref{dnls}) blows up in time $t=T$; then
$$
\lim_{t\uparrow T} \|u(t)\|_{L^6}=\infty.
$$
\end{remark}

\begin{remark}
It is unclear whether (\ref{eq:assum}) holds, and the question is an open problem. 
\end{remark}

\begin{notation}
Throughout the paper, we use letters $c,~C,~c(*),~C(*)$ to denote various constants.
\end{notation}

The paper is organized as follows.
Section \ref{sec:p} describes the analysis of the gauge group of transformations associated to (\ref{dnls}).
In Section \ref{sec:H1estimate}, we derive the a priori estimate for solutions in $\dot{H}^1(\mathbb{R})$ space by using a contradiction argument.
In Subsection \ref{sec:concentration}, we perform qualitative analysis to obtain some concentration phenomena; the blow-up solutions behave like a standing wave solution.
In Subsection \ref{sec:proofthm}, we provide the proof of Theorem \ref{thm:main}. 
Section \ref{sec:exam} gives the proof of Proposition \ref{prop:exam}.
Section \ref{sec:P2Q} is devoted to the direct calculation of $E_2(q_{c^2/4,c})=P_2(q_{c^2/4,c})=P_3(q_{c^2/4,c})=0$ as an appendix.

\section{Preliminary}\label{sec:p}

We begin by defining the gauge transformation associated to (\ref{dnls}).
If we assume $u$ to be a solution to the Cauchy problem (\ref{dnls}), then
$$
{\mathcal G}_{\nu}u(t,x)=u(t,x)e^{-i\frac{\nu}{2}\int_{-\infty}^x|u(t,y)|^2\,dy},\quad \nu\in\mathbb{R}
$$
is formally a solution of the following Cauchy problem;
\begin{eqnarray}\label{dnls-g}
\left\{
\begin{array}{l}
i\partial_tv+\partial_x^2v=i(2-\nu)|v|^2\partial_xv+i(1-\nu)v^2\partial_x\overline{v}+\frac{\nu(1-\nu)}{4}|v|^4v,\\
v(0,x)=v_{0,\nu}(x), 
\end{array}
\right.
\end{eqnarray}
where $v_{0,\nu}={\mathcal G}_{\nu}(u_0)$.
Using the gauge transform, the conservation laws in (\ref{eq:mass}), (\ref{eq:H1energy}), (\ref{eq:H1momentum}), (\ref{eq:H2energy}) and (\ref{eq:H2momentum}) reduce to the following, the $L^2$-norm
\begin{eqnarray}\label{eq:gmass}
{\mathcal M}(v)=M(u),
\end{eqnarray}
the energy
\begin{eqnarray}\label{eq:gH1energy}
{\mathcal E}_{1,\nu}(v) & = & \int_{-\infty}^{\infty}\left(|\partial_x v|^2-\left(\frac32-\nu\right)|v|^2\Im\left(\overline{v}\partial_xv\right)+\frac{(2-\nu)(1-\nu)}{4}|v|^6\right)\,dx\nonumber\\
& = & {\mathcal E}_{1,\nu}(v_{0,\nu})=E(u_0),
\end{eqnarray}
the momentum
\begin{eqnarray}\label{eq:gH1momentum}
{\mathcal P}_{1,\nu}(v)=\int_{-\infty}^{\infty}\left(\Im\left(\overline{v}\partial_x v\right)-\frac{1-\nu}{2}|v|^4\right)\,dx={\mathcal P}_{1,\nu}(v_{0,\nu})=P_1(u_0),
\end{eqnarray}
the $H^2$-energy
\begin{eqnarray}\label{eq:gH2energy}
{\mathcal E}_{2,\nu}(v) & = & E_2(v)+\frac{\nu}{16}\left(\nu^3-10\nu^2+30\nu-35\right)\|v\|_{L^{10}}^{10}
+\nu(4\nu-15)\||v|^2\partial_xv\|_{L^2}^2\\
& & 
+\frac52\nu(\nu-3)\Re\langle(\partial_xv)^2,|v|^2v^2\rangle
+3\nu\Im\langle \partial_x^2v,|v|^2\partial_xv\rangle+\nu\Im\langle \partial_x^2v,v^2\overline{\partial_xv}\rangle\nonumber\\
& & 
+\frac{\nu}{4}(2\nu^2-15\nu+30)\Im\langle\partial_xv,|v|^6v\rangle\nonumber\\
& = & {\mathcal E}_{2,\nu}(v_{0,\nu}) \nonumber\\
& = & E_2(u_0)\nonumber
\end{eqnarray}
and the $H^2$-momentum
\begin{eqnarray}\label{eq:gH2momentum}
{\mathcal P}_{2,\nu}(v) & = & P_2(v)+\frac{\nu}{16}(\nu^2-6\nu+10)\|v\|_{L^8}^8+\frac54\nu\|v\partial_xv\|_{L^2}^2\\
& & +\frac{\nu}{2}\Re\langle(\partial_xv)^2,v^2\rangle+\frac38\nu(\nu-4)\Im\langle \partial_xv,|v|^4v\rangle\nonumber\\
& = & {\mathcal P}_{2,\nu}(v_{0,\nu})\nonumber\\
& = & P_2(u_0).\nonumber
\end{eqnarray}
In particular, the equation of (\ref{dnls-g}) with parameter $\nu=3/2$ is
\begin{eqnarray}\label{dnls-gp}
\left\{
\begin{array}{l}
i\partial_tv+\partial_x^2v=\frac{i}{2}\left(|v|^2\partial_xv-v^2\partial_x\overline{v}\right)-\frac{3}{16}|v|^4v,\\
v(0,x)=v_0, 
\end{array}
\right.
\end{eqnarray}
where $v_0=v_{0,3/2}$.

It is easy to see that the map ${\mathcal G}_{\nu}$ is a bijection in $H^s$ for $s\ge 1/2$ (see \cite{t1}).
Let us consider $v_0\in L^2$ such that the Cauchy problem (\ref{dnls-gp}) blows up at time $t=T$.

\section{$\dot{H}^1$-norm estimate for solutions to the gauged equivalent equation}\label{sec:H1estimate}

First, we recall that because of the time reflection invariance, it suffices to consider non-negative time.

Let $v$ be the $H^1$-solution to (\ref{dnls-gp}).
If $\sup_{0\le t< T}\|v(t)\|_{L^6}\le C(T)$ for any $T\in(0,\infty]$, then we have from the energy conservation law (\ref{eq:gH1energy}) that $\sup_{0\le t<T}\|\partial_xv(t)\|_{L^2}\le C(T)$ for any $T>0$.
In addition, in \cite{fhi} the global $\dot{H}^1$ bound of solutions was obtained for ${\mathcal P}_{1,3/2}(v_0)<0$.
Then without loss of generality, we assume that ${\mathcal P}_{1,3/2}(v_0)\ge 0$ and $\limsup_{t\to T-}\|v(t)\|_{L^6}=\infty$ for some time $T\in(0,\infty]$.

Given solutions $v={\mathcal G}_{3/2}u$ of (\ref{dnls-gp}), for $\alpha\in\mathbb{R}$ we write $w(t,x)=v(t,x)e^{i\alpha x}$, which satisfies the following equation;
$$
i\partial_tw+\partial_x^2w-2i\alpha\partial_xw-\alpha^2w=\frac{i}{2}\left(|w|^2\partial_xw-w^2\partial_x\overline{w}\right)+\alpha|w|^2w-\frac{3}{16}|w|^4w.
$$
In this setting, the conservation laws ${\mathcal P}_{1/3/2}$ and ${\mathcal E}_{1,3/2}$ become
\begin{eqnarray*}
{\cal P}_{1,3/2}(v)=\Im\int_{-\infty}^{\infty}\overline{w}\partial_xw\,dx+\frac{1}{4}\|w\|_{L^4}^4-\alpha\|w\|_{L^2}^2
\end{eqnarray*}
and
\begin{eqnarray*}
{\cal E}_{1/3/2}(v)=\|\partial_xw\|_{L^2}^2-\frac{1}{16}\|w\|_{L^6}^6-2\alpha\Im\int_{-\infty}^{\infty}\overline{w}\partial_xw\,dx+\alpha^2\|w\|_{L^2}^2,
\end{eqnarray*}
respectively, that imply
\begin{eqnarray}\label{pe}
{\cal E}_{1,3/2}(v)=\|\partial_xw\|_{L^2}^2-\frac{1}{16}\|w\|_{L^6}^6-2\alpha\left({\cal P}_{1,3/2}(v)-\frac14\|w\|_{L^4}^4\right)-\alpha^2\|w\|_{L^2}^2.
\end{eqnarray}

We use some ideas from \cite{a1}.
The following sharp Gagliardo-Nirenberg inequality was obtained (see also \cite{gw}). 
\begin{lemma}\label{lem:GN}
\begin{eqnarray*}
\|\phi\|_{L^6}\le C_{GN}  \|\partial_x\phi\|_{L^2}^{1/9}\|\phi\|_{L^4}^{8/9},
\end{eqnarray*}
where $C_{GN}=3^{1/6}\cdot(2\pi)^{-1/9}$.
\end{lemma}

If $\alpha>0$, we apply Lemma \ref{lem:GN} to obtain the following
\begin{eqnarray*}
\frac{\|v\|_{L^4}^4}{4}+\frac{C_{GN}^{-18}\|v\|_{L^6}^{18}}{2\alpha\|v\|_{L^4}^{16}}\le \frac{{\cal E}_{1,3/2}(v)}{2\alpha}+{\cal P}_{1,3/2}(v)+\frac{\alpha\|v\|_{L^2}^2}{2}+\frac{\|v\|_{L^6}^6}{32\alpha}.
\end{eqnarray*}
By scaling the parameter $\alpha$ as $\alpha\|v\|_{L^6}^3$, we have
\begin{eqnarray}\label{eq:gw22}
\left(\frac{\|v\|_{L^4}^4}{4\|v\|_{L^6}^3}+\frac{C_{GN}^{-18}\|v\|_{L^6}^{15}}{2\alpha\|v\|_{L^4}^{16}}-\frac{\|v\|_{L^2}^2}{2}\alpha-\frac{1}{32\alpha}\right)\|v\|_{L^6}^3\le  \frac{{\cal E}_{1,3/2}(v)}{2\alpha\|v\|_{L^6}^3}+{\cal P}_{1,3/2}(v).\end{eqnarray}
Now we put $X=\|v\|_{L^4}^4/\|v\|_{L^6}^3$, and choose $\alpha=\alpha_0=2^{-5/2}\cdot 3^{-1/2}\cdot \pi^{-1/2}$.
Then the left-hand side of (\ref{eq:gw22}) is
$$
\left(\frac{X}{4}+\frac{C_{GN}^{-18}}{2\alpha_0X^4}-\frac{\|v\|_{L^2}^2}{2}\alpha_0-\frac{1}{32\alpha_0}\right)\|v(t)\|_{L^6}^3.
$$
Observe that upon using the assumption $\|v\|_{L^2}=\sqrt{4\pi}$ in Theorem \ref{thm:main}, the function
$$
f(x)=\frac{x}{4}+\frac{C_{GN}^{-18}}{2\alpha_0x^4}-\frac{\|v\|_{L^2}^2}{2}\alpha_0-\frac{1}{32\alpha_0},\quad x>0
$$
has a minimum value zero at
$$
x=X_0=\left(\frac{8C_{GN}^{-18}}{\alpha_0}\right)^{1/5}=2^{3/2}\cdot 3^{-1/2}\cdot\pi^{1/2},
$$
namely $f(X_0)=0$.
Note that
$$
\frac{x}{4}+\frac{C_{GN}^{-18}}{2\alpha_0x^4}-2\pi\alpha_0-\frac{1}{32\alpha_0}=0,
$$
is possible only if $x=X_0$.
We reapply the above estimate to (\ref{eq:gw22}).
Consequently, it follows that  
\begin{eqnarray}\label{eq:bound}
0\le \left(\frac{X}{4}+\frac{C_{GN}^{-18}}{2\alpha_0X^4}-\frac{\|v\|_{L^2}^2}{2}\alpha_0-\frac{1}{32\alpha_0}\right)\|v\|_{L^6}^3\le \frac{{\cal E}_{1,3/2}(v)}{2\alpha_0\|v\|_{L^6}^3}+{\cal P}_{1,3/2}(v).
\end{eqnarray}

If there exists a time $T\in(0,\infty]$ such that 
$$
\liminf_{t\to T-}\left|\frac{\|v(t)\|_{L^4}^4}{\|v(t)\|_{L^6}^3}-X_0\right|>0,
$$ 
then we have
$$
\liminf_{t\to T-}\left(\frac{X}{4}+\frac{C_{GN}^{-18}}{2\alpha_0X^4}-\frac{\|v\|_{L^2}^2}{2}\alpha_0-\frac{1}{32\alpha_0}\right)>0,
$$
which leads $\limsup_{t\to T-}\|v(t)\|_{L^6}<\infty$.
It is not deemed acceptable. 
 
As we have seen above, blow-up solutions $v$ should satisfy
$$
\liminf_{t\to T-}\left|\frac{\|v(t)\|_{L^4}^4}{\|v(t)\|_{L^6}^3}-X_0\right|=0.
$$
Then there exists a sequence $\{t_n\}$ such that
$$
t_n\le t_{n+1}<T,\quad \lim_{n\to\infty}t_n=T,\quad \lim_{n\to\infty}\frac{\|v(t_n)\|_{L^4}^4}{\|v(t_n)\|_{L^6}^3}=X_0,
$$
namely
$$
\frac{\|v(t_n)\|_{L^4}^4}{\|v(t_n)\|_{L^6}^3}=X_0+o_n(1),
$$
where $\lim_{n\to\infty}o_n(1)=0$.
Using the conservation law ${\mathcal E}_{1,3/2}$,  one can get $\lim_{n\to\infty}\|\partial_xv(t_n)\|_{L^2}=\infty$.
Then it follows that
$$
\|v(t_n)\|_{L^4}^4=X_0\|v(t_n)\|_{L^6}^3+o_n(1)\|v(t_n)\|_{L^6}^3.
$$

For $c>0$ fixed, which is the presumed constant, see Remark \ref{rem:const} below, take
$$
\lambda_n=\sqrt{\frac{3\cdot 2^3\cdot c^2\cdot \pi}{\|v(t_n)\|_{L^6}^6}},
$$
which satisfies $\lim_{n\to \infty}\lambda_n=0$.
We rescale the solution $v$ by
$$
w_n(x)=\lambda_n^{1/2}v(t_n,\lambda_nx)e^{i\beta x},
$$
where $\beta\in\mathbb{R}$ is chosen later.
Hence, we know
$$
\|w_n\|_{L^4}^4=\lambda_n\|v(t_n)\|_{L^4}^4,\quad \|w_n\|_{L^6}^3=\lambda_n \|v(t_n)\|_{L^6}^3
$$
and
$$
\|\partial_xw_n\|_{L^2}^2=\lambda_n^2\|\partial_xv(t_n)\|_{L^2}^2+2\beta\lambda_n\Im\int_{\mathbb{R}}\overline{v(t_n,x)}\partial_xv(t_n,x)\,dx+\beta^2M(v_0).
$$
In particular we have
$$
\|w_n\|_{L^6}^6=3\cdot 2^3\cdot c^2\cdot \pi=\|Q_c\|_{L^6}^6,
$$
$$
\|w_n\|_{L^4}^4=\lambda_n\left(X_0\|v(t_n)\|_{L^6}^3+o_n(1)\|v(t_n)\|_{L^6}^3\right)\to X_0\cdot 3^{\frac12}\cdot 2^{\frac32}\cdot c \cdot \pi^{\frac12}=8c\pi=\|Q_c\|_{L^4}^4,
$$
as $n\to \infty$.
Moreover choosing $\beta=c/2$, one gets
\begin{eqnarray*}
\|\partial_xw_n\|_{L^2}^2& = &\lambda_n^2\left({\cal E}_{1,3/2}(v_0)+\frac{1}{16}\|v(t_n)\|_{L^6}^6\right)+2\beta\lambda_n\left({\cal P}_{1,3/2}(v_0)-\frac{1}{4}\|v(t_n)\|_{L^4}^4\right)+4\pi\beta^2\\
& \to  & \frac{3\cdot 2^3\cdot c^2 \cdot \pi}{2^4}-\frac{2\cdot \beta\cdot 3^{\frac12}\cdot 2^{\frac32}\cdot c\cdot \pi^{\frac12}\cdot 2^{\frac32}\cdot 3^{-\frac12}\cdot \pi^{\frac12}}{2^2}+4\pi\beta^2\\
& = & \frac{\pi}{2}\left(8\left(\beta-\frac{c}{2}\right)^2 +c^2\right)\\
& = & \frac{c^2\pi}{2}=\|\partial_xQ_c\|_{L^2}^2,
\end{eqnarray*}
as $n$ goes to infinity.

As conclusion, putting
$$
w_n(x)=\lambda_n^{1/2}v(t_n,\lambda_nx)e^{icx/2},
$$
we obtain
$$
\lim_{n\to\infty}\|w_n\|_{L^4}^4=\|Q_c\|_{L^4}^4,\quad \|w_n\|_{L^6}^6=\|Q_c\|_{L^6}^6,\quad \lim_{n\to\infty}\|\partial_xw_n\|_{L^2}^2=\|\partial_xQ_c\|_{L^2}^2.
$$

\begin{remark}
We comment on the profile of $w_n$.
Recall
$$
v(t_n,x)=\frac{1}{\lambda_n^{1/2}}e^{-i\frac{c}{2\lambda_n}x}w_n\left(\frac{x}{\lambda_n}\right).
$$
It follows 
\begin{eqnarray*}
{\cal P}_{1,3/2}(v(t_n)) & =& \frac{1}{\lambda_n^2}\Im\int_{\mathbb{R}}\overline{w_n\left(\frac{x}{\lambda_n}\right)}\left\{\partial_xw_n\left(\frac{x}{\lambda_n}\right)-i\frac{c}{2}w_n\left(\frac{x}{\lambda_n}\right)\right\}\,dx+\frac{1}{4\lambda_n}\|w_n\|_{L^4}^4\\
& = & \frac{1}{\lambda_n}\left({\mathcal P}_{1,3/2}(w_n)-\frac{c}{2}\|w_n\|_{L^2}^2\right)
\end{eqnarray*}
and then
\begin{eqnarray}\label{eq:scvw}
{\cal E}_{1,3/2}(v(t_n))  & = & \frac{1}{\lambda_n^3}\left\|\partial_xw_n\left(\frac{\cdot}{\lambda_n}\right)-i\frac{c}{2}w_n\left(\frac{\cdot}{\lambda_n}\right)\right\|_{L^2}^2-\frac{1}{16\lambda_n^2}\|w_n\|_{L^6}^6\\
& = & \frac{1}{\lambda_n^2}\left(\|\partial_xw_n\|_{L^2}^2+\frac{c^2}{4}\|w_n\|_{L^2}^2-c\Im\int_{\mathbb{R}}\overline{w_n(x)}\partial_xw_n(x)\,dx-\frac{1}{16}\|w_n\|_{L^6}^6\right)\nonumber\\
& = & \frac{1}{\lambda_n^2}\left({\mathcal E}_{1,3/2}(w_n)
+\frac{c}{4}\|w_n\|_{L^4}^4-\frac{c^2}{4}\|w_n\|_{L^2}^2-c{\cal P}_{1,3/2}(v(t_n))\lambda_n\right).\nonumber
\end{eqnarray}
It is important to notice that since $\|w_n\|_{L^6}=\|Q_c\|_{L^6}$ together with (\ref{eq:scvw}), we have
\begin{eqnarray*}
\frac{5}{2}c^2\pi
& = & \lim_{n\to\infty}\left(\|\partial_xw_n\|_{L^2}^2+\frac{c}{4}\|w_n\|_{L^4}^4\right)\\
& = & \lim_{n\to\infty}\left(\frac{1}{16}\|w_n\|_{L^6}^6+\frac{c^2}{4}\|w_n\|_{L^2}^2+c{\cal P}_{1,3/2}(v(t_n))\lambda_n+{\cal E}_{1,3/2}(v(t_n))\lambda_n^2\right)\\
& = & \lim_{n\to\infty}\left(\frac{5}{2}c^2\pi+c{\cal P}_{1,3/2}(v_0)\lambda_n+{\cal E}_{1,3/2}(v_0)\lambda_n^2\right).
\end{eqnarray*}
Observe that using (\ref{eq:minmizer}), we get that there holds either (i) ${\cal P}_{1,3/2}(v_0)> 0$ or (ii) ${\cal P}_{1,3/2}(v_0)=0<{\cal E}_{1,3/2}(v_0)$; in the other cases, the failure would lead to the global well-posedness.
\end{remark}

\section{Proof of Theorem \ref{thm:main}}

\subsection{Concentration phenomena}\label{sec:concentration}

To obtain a concentration phenomenon and a singularity profile of solutions, we need the following lemma.

\begin{lemma}[Corollary 3.1 \cite{a1}]\label{lem:Au}
\begin{eqnarray*}
\min\left\{\frac12\|\partial_xu\|_{L^2}^2+\frac14\|u\|_{L^4}^4\mid\|u\|_{L^6}=1\right\}=\frac{5}{18}\lambda,
\end{eqnarray*}
where the minimum is attainted by, without gauge invariance and dilation,
$$
u^*(x)=\sqrt{\frac{2}{x^2+\frac{4}{3}\lambda}},
$$
for $\lambda=3\cdot (3\pi)^{2/5}/4$ satisfying $\|u^*\|_{L^6}=1$, which satisfies
$$
-u''+u^3-\lambda u^5=0.
$$
\end{lemma}

Now if
$$
\mu=2\cdot(3\pi)^{\frac15}\cdot\nu^2,\quad
\nu^2=\frac{(3\pi)^{1/10}}{2\cdot 3^{1/2}\cdot \pi^{1/2}\cdot c},
$$
one has
$$
\|w_n\|_{L^6}^6=\frac{\mu}{\nu^6},\quad
\left(\frac{\nu}{\mu}\right)^2=\frac{c}{2},\quad
\frac34\cdot(3\pi)^{2/5}\cdot\frac{\nu^4}{\mu^2}=\frac{3}{16}.
$$
We use these to Lemma \ref{lem:Au}.
Using $u(x)=\nu v(\mu x)$, it follows
\begin{eqnarray*}
\min\left\{\frac{\nu^2\mu}{2}\|\partial_xv\|_{L^2}^2+\frac{\nu^4\mu^{-1}}{4}\|v\|_{L^4}^4\mid\|v\|_{L^6}^6=\frac{\mu}{\nu^6}\right\}=\frac{5}{18}\lambda.
\end{eqnarray*}
Hence
\begin{eqnarray}\label{eq:minmizer}
\min\left\{\|\partial_xv\|_{L^2}^2+\frac{c}{4}\|v\|_{L^4}^4\mid\|v\|_{L^6}^6=\|Q_c\|_{L^6}^6\right\}=\frac{5}{2}c^2\pi,
\end{eqnarray}
where the minimum is attainted by 
$$
v(x)=\frac{1}{\nu}u^*\left(\frac{x}{\mu}\right),
$$
which satisfies the following equation
$$
-v''+\frac{c}{2}v^3-\frac{3}{16} v^5=0.
$$
In particular, the property in (\ref{eq:minmizer}) allows us to establish that for all $\varepsilon>0$, there exists a natural number $N\in\mathbb{N}$ such that for $n\ge N$, we have
$$
\frac{5}{2}c^2\pi<\|\partial_xw_n\|_{L^2}^2+\frac{c}{4}\|w_n\|_{L^4}^4<\frac{5}{2}c^2\pi+\varepsilon.
$$

Recall the following profile decomposition lemma (see e.g. \cite{g,hk,km,k}).

\begin{lemma}[The profile decomposition (e.g. \cite{g,hk,km,k})]\label{lem:profile}
There exist $\{V^j\},~\{x_n^j\}$ such that 
$$
w_n(x)=\sum_{j=1}^LV^j(x-x_n^j)+R_n^L(x),
$$
where for $j\ne k$
\begin{eqnarray}\label{eq:diff}
\lim_{n\to\infty}|x_n^j-x_n^k|=\infty,
\end{eqnarray}
\begin{eqnarray}\label{eq:r-mass}
\|Q_c\|_{L^2}^2=\|w_n\|_{L^2}^2=\sum_{j=1}^L\|V^j\|_{L^2}^2+\|R_n^L\|_{L^2}^2+o_n(1),
\end{eqnarray}
\begin{eqnarray}\label{eq:r-kin}
\|\partial_xw_n\|_{L^2}^2=\sum_{j=1}^L\|\partial_xV^j\|_{L^2}^2+\|\partial_xR_n^L\|_{L^2}^2+o_n(1)=\|\partial_xQ_c\|_{L^2}^2+o_n(1),
\end{eqnarray}
\begin{eqnarray}\label{eq:r-energy}
{\mathcal E}_{1,3/2}(w_n)=\sum_{j=1}^L{\mathcal E}_{1,3/2}(V^j)+{\mathcal E}_{1,3/2}(R_n^L)+o_n(1)={\mathcal E}_{1,3/2}(Q_c)+o_n(1),
\end{eqnarray}
and
\begin{eqnarray}\label{eq:r-L6}
\limsup_{n,L\to\infty}\|R_n^L\|_{L^p}=0
\end{eqnarray}
for every $p\in(2,\infty)$.
\end{lemma}

Note that in view of (\ref{eq:diff}), (\ref{eq:r-mass}), (\ref{eq:r-kin}) and (\ref{eq:r-L6}), this lemma implies
\begin{eqnarray}\label{eq:in}
\|w_n\|_{L^p}^p=\sum_{j=1}^L\|V^j\|_{L^p}^p+\|R_n^L\|_{L^p}^p+o_n(1)=\sum_{j=1}^L\|V^j\|_{L^p}^p+o_{n,L}(1),
\end{eqnarray}
for $p=4$ and $6$, where $\lim_{n,L\to\infty}o_{n,L}(1)=0$.
Now let $K_c$ be the functional given by
$$
K_c(f)=\|\partial_xf\|_{L^2}^2+\frac{c}{4}\|f\|_{L^4}^4.
$$
Since $\|V^j\|_{L^6}\le \|w_n\|_{L^6}$, we have
\begin{eqnarray}\label{eq:K}
\frac{5}{2}c^2\pi+\varepsilon & > & K_c(w_n)\\
& = & \sum_{j=1}^LK_c(V_j)+K_c(R_n^L)+o_n(1)\nonumber\\
& \ge & \sum_{j=1}^LK_{c_j}(V^j)+o_n(1)\nonumber\\
& \ge & \sum_{j=1}^L\frac{5}{2}c_j^2\pi+o_n(1),\nonumber
\end{eqnarray}
where $c_j\in[0,c]$ is defined by $\|V^j\|_{L^6}^6=3\cdot 2^3\cdot c_j^2\cdot \pi$.
Also, from (\ref{eq:in}), the right-hand side is
$$
\sum_{j=1}^L\frac{5}{2}c_j^2\pi+o_n(1)=\|w_n\|_{L^6}^6+o_{n,L}(1)=\frac{5}{2}c^2\pi+o_{n,L}(1),
$$
which implies that
$$
K_{c_j}(V^j)=\frac{5}{2}c_j^2\pi
$$
for all $1\le j\le L$.
Since $\|V^j\|_{L^6}^6=3\cdot 2^3\cdot c_j^2\cdot \pi$, by Lemma \ref{lem:Au}, we obtain that $V^j(x)=Q_{c_j}(x-x_j)e^{i\alpha_j}$ for some $(x_j,\alpha_j)\in\mathbb{R}^2$.
Furthermore, by (\ref{eq:r-mass}), we obtain
$$
4\pi=\sum_{j=1}^L\|Q_{c_j}\|_{L^2}^2+\|R_{n}^L\|_{L^2}^2+o_{n}(1),
$$
for which one has $L=1$, because of $\|Q_{c_j}\|_{L^2}=\sqrt{4\pi}$.

\subsection{Proof of Theorem \ref{thm:main}}\label{sec:proofthm}

Collecting the information above, we describe $w_n$ as
$$
w_n(x)=e^{i\gamma}\left(Q_c(x-x_n)+R_n(x)\right).
$$
for some $\gamma\in\mathbb{R}$ and $x_n\in\mathbb{R}$. 
As a consequence, we have
\begin{eqnarray}\label{eq:profile}
v(t_n,x)=e^{-i\frac{c}{2\lambda_n}x+i\gamma}\frac{1}{\lambda_n^{1/2}}\left(Q_c\left(\frac{x}{\lambda_n}-x_n\right)+R_n\left(\frac{x}{\lambda_n}\right)\right),
\end{eqnarray}
where $R_n$ satisfies
$$
\|R_n\|_{L^2}^2+2\Re\langle Q_c(\cdot-x_n),R_n\rangle=0.
$$
Moreover
$$
\left\|v(t_n,x)-e^{-i\frac{c}{2\lambda_n}x+i\gamma}\frac{1}{\lambda_n^{1/2}}R_n\left(\frac{x}{\lambda_n}\right)\right\|_{L^p}=\|Q_c\|_{L^p}=\|v(t_n)\|_{L^p},
$$
for $p=2$ and $p=6$.

Denote $v_n=v(t_n,\cdot)$, which completes the proof of Theorem \ref{thm:main}.
\qed
 
\begin{remark}\label{re:K}
As a consequence of the above calculation in (\ref{eq:K}), we get
$$
\lim_{n\to\infty}K_c(R_n)=0.
$$
\end{remark}
\begin{remark}\label{rem:const}
Notice that if we replace $\lambda_n$ by $\lambda_n/c$ (absorbed into $\lambda_n$), then we ignore the irrelevant constant $c$, and (\ref{eq:profile}) can be also written as the form
$$
v(t_n,x)=e^{-i\frac{x}{2\lambda_n}+i\gamma}\frac{1}{\lambda_n^{1/2}}\left(Q_1\left(\frac{x}{\lambda_n}-x_n\right)+R_n\left(\frac{x}{\lambda_n}\right)\right).
$$
\end{remark}

\section{Proof of Proposition \ref{prop:exam}}\label{sec:exam}

In this section, we prove Proposition \ref{prop:exam} by giving an example.

Let $q_0=q_{1/4,1}\in H^3(\mathbb{R})$ and define
$$
\varphi_a(x)=\frac{\sqrt{x^2+1}}{2(x^4+1)}e^{-\frac{i}{2}x+3i\left(\arctan x+\frac{\pi}{2}\right)-iax}\in H^3(\mathbb{R}),\quad a\in \mathbb{R}.
$$
Notice that
$$
\|\varphi_a\|_{L^2}^2=\frac14\int_{\mathbb{R}}\frac{1+x^2}{(1+x^4)^2}\,dx\in(0,\infty)
$$
is independent of $a\in\mathbb{R}$.

A simple computation shows the following lemma.

\begin{lemma}\label{lem:com}
For $a\in\mathbb{R}$ and $b>0$,
$$
\int_{\mathbb{R}}\frac{\cos ax}{1+x^4}\,dx=\pi e^{-\frac{a}{\sqrt{2}}}\sin\left(\frac{a}{\sqrt{2}}+\frac{\pi}{4}\right),
$$
$$
\int_{\mathbb{R}}\frac{\cos ax}{(x^2+1)(x^4+1)}\,dx=\frac{\pi}{2}\left(e^{-a}+\sqrt{2}e^{-\frac{a}{\sqrt{2}}}\sin\frac{a}{\sqrt{2}}\right),
$$
$$
\int_{\mathbb{R}}\frac{x\sin ax}{(x^2+1)(x^4+1)}\,dx=\frac{\pi}{2}\left(e^{-a}-e^{-\frac{a}{\sqrt{2}}}\cos\frac{a}{\sqrt{2}}+e^{-\frac{a}{\sqrt{2}}}\sin\frac{a}{\sqrt{2}}\right).
$$
and
$$
\int_{\mathbb{R}}\frac{\cos( bx-3\arctan x)\sqrt{x^2+1}}{x^4+1}\,dx=2\pi\left(-2e^{-b}+\frac14\cos\frac{b}{\sqrt{2}}e^{-\frac{b}{\sqrt{2}}}+\frac14\left(1-\frac{2}{\sqrt{2}}\right)\sin\frac{b}{\sqrt{2}}e^{-\frac{b}{\sqrt{2}}}\right).
$$
\end{lemma}

Consider
$$
u_0(x)=q_0(x)+\eta\varphi_a(x).
$$
where $a,\eta\in\mathbb{R}$ will be chosen later. 
We illustrate that $u_0$ satisfies (\ref{eq:example}) for some $a$ and $\eta$. 

Since $\|q_0\|_{L^2}=\sqrt{4\pi}$, we easily see that $\|u_0\|_{L^2}=\sqrt{4\pi}$, if and only if
\begin{eqnarray}\label{eq:con11}
\eta\|\varphi_a\|_{L^2}^2+2\Re\langle q,\varphi_a\rangle=0.
\end{eqnarray}
Using Lemma \ref{lem:com}, one obtains that (\ref{eq:con11}) is equivalent to the following:
\begin{eqnarray}\label{eq:con1}
\eta=-\frac{2\pi}{\|\varphi_a\|_{L^2}^2}e^{-\frac{a}{\sqrt{2}}}\sin\left(\frac{a}{\sqrt{2}}+\frac{\pi}{4}\right).
\end{eqnarray}
So, we pick up $a>0$ such that
\begin{eqnarray}\label{eq:con2}
\frac{a}{\sqrt{2}}+\frac{\pi}{4}=2k\pi-\varepsilon_0
\end{eqnarray}
for $0<\varepsilon_0\ll 1$ and $k\in\mathbb{Z}$, by demanding
\begin{eqnarray}\label{eq:con3}
0<\eta\lesssim e^{-\frac{a}{\sqrt{2}}}\varepsilon_0\ll 1.
\end{eqnarray}

Observe that
\begin{eqnarray}\label{eq:con21}
\partial_xq_0=-\frac{i}{2}q_0+\frac34iQ_1^2q_0-\frac{x}{x^2+1}q_0
\end{eqnarray}
and
\begin{eqnarray}\label{eq:con22}
\partial_x\varphi_a=-\frac{i}{2}\varphi_a+\frac34iQ_1^2\varphi_a-ia\varphi_a-\frac{4x^3}{x^4+1}\varphi_a+\frac{x}{\sqrt{x^2+1}}\varphi_a.
\end{eqnarray}
Using $P_1(q_0)=0$, it follows that
$$
P_1(u_0)=2\eta\Im\langle \partial_xq_0,\varphi_a\rangle+\eta^2\Im\langle \partial_x\varphi_a,\varphi_a\rangle-\frac12\left(\|q_0+\eta\varphi_a\|_{L^4}^4-\|q_0\|_{L^4}^4\right).
$$
From (\ref{eq:con11}), (\ref{eq:con21}), (\ref{eq:con22}) and (\ref{eq:con3}), we have that $0<\eta\ll 1$ and
$$
P_1(u_0)=-2\eta\left(\frac14\Re\langle Q_1^2q_0,\varphi_a\rangle+\Im\left\langle \frac{x}{x^2+1}q_0,\varphi_a\right\rangle\right)+O(\eta^2).
$$
By Lemma \ref{lem:com}, one deduces that
\begin{eqnarray*}
P_1(u_0) & = & -\eta\pi\left(2e^{-a}+e^{-\frac{a}{\sqrt{2}}}\left((\sqrt{2}+1)\sin\frac{a}{\sqrt{2}}-\cos\frac{a}{\sqrt{2}}\right)\right)+O(\eta^2)\nonumber\\
& = & -\eta\pi\left(2e^{-a}+e^{-\frac{a}{\sqrt{2}}}\sqrt{4+2\sqrt{2}}\sin\left(\frac{a}{\sqrt{2}}-\theta_0\right)\right)+O(\eta^2),
\end{eqnarray*}
where $\theta_0\in(0,\pi/2)$ such that
$$
\cos\theta_0=\frac{\sqrt{2}+1}{\sqrt{4+2\sqrt{2}}},\quad
\sin\theta_0=\frac{1}{\sqrt{4+2\sqrt{2}}}.
$$
For large $k\in\mathbb{N}$ in (\ref{eq:con2}), we have $e^{-a}\ll e^{-a/\sqrt{2}}$ and
$$
\sin\left(\frac{a}{\sqrt{2}}-\theta_0\right)=-\sin\left(\theta_0+\frac{\pi}{4}-\varepsilon_0\right)\lesssim -1.
$$
Therefore
$$
P_1(u_0)\gtrsim \eta e^{-\frac{a}{\sqrt{2}}}+O(\eta^2)>0.
$$
In particular, since $M_0(q_0)=0$ and by Lemma \ref{lem:com}, one gets $M_0(u_0)\ne 0$, if $\eta>0$ is small.

This completes the proof of Proposition \ref{prop:exam}.
\qed 

\section{Appendix}\label{sec:P2Q}

In this section we calculate $E_2(q_{c^2/4_c})=P_2(q_{c^2/4,c})=P_3(q_{c^2/4,c})=0$.

We begin with summarizing a few computations;
\begin{eqnarray}\label{eq:com}
& & \|Q_c\|_{L^2}^2=2^2\pi,~\|Q_c\|_{L^4}^4=2^3c\pi,~\|Q_c\|_{L^6}^6=3\cdot 2^3c^2\pi,~\|Q_c\|_{L^8}^8=5\cdot 2^4c^3\pi,\\
& & \|Q_c\|_{L^{10}}^{10}=7\cdot 5\cdot 2^3c^4\pi,~\|Q_c\|_{L^{12}}^{12}=9\cdot 7\cdot 2^4c^5\pi,\nonumber\\
& & \|\partial_xQ_c\|_{L^2}^2=\frac12c^2\pi,~\|Q_c\partial_xQ_c\|_{L^2}^2=c^3\pi,~\|Q_c^2\partial_xQ_c\|_{L^2}^2=\frac{5}{2}c^4\pi,~\|Q_c^3\partial_xQ_c\|_{L^2}^2=7c^5\pi,\nonumber\\
& & \|\partial_xQ_c\|_{L^4}^4=\frac{3}{2^4}c^5\pi.\nonumber
\end{eqnarray}

For simplicity, we shall only consider the case  $c=1$.

Recall the recurrence formula obtained in \cite{tsm1} (also in \cite{kn}).
Let $\{Z^{(n)}\}$ be determined by 
\begin{eqnarray}\label{eq:rec}
Z^{(n+1)}=\sum_{k=1}^{n-1}Z^{(k)}Z^{(n-k)}+iqrZ^{(n)}+q\frac{\partial}{\partial x}\left(\frac1qZ^{(n)}\right),\quad n\in\mathbb{N},
\end{eqnarray}
where $Z^{(0)}=qr$ and
\begin{eqnarray}\label{eq:z1}
Z^{(1)}=-\frac14q^2r^2+\frac{i}{2}qr_x.
\end{eqnarray}
In \cite{kn,tsm1}, it was proved that if $q=q(t,x)$ is a solution to (\ref{dnls}) and $r=\overline{q}$, every $\int_{\mathbb{R}}Z^{(i)}\,dx$ is a conserved functional.
In particular, we have that for a solution $q=q(t,x)$ to (\ref{dnls}) and $r=\overline{q}$, 
$$
M(q)=\int_{\mathbb{R}}Z^{(0)}\,dx,\quad
P_n(q)=2
\int_{\mathbb{R}}Z^{(2n-1)}\,dx,\quad
E_n(q)=-2
\int_{\mathbb{R}}Z^{(2n)}\,dx,\quad n\in\mathbb{N}.
$$
The density $Z^{(n)}~(n=2,3,4)$ are computed in \cite{tsm1} as follows:  
\begin{eqnarray}\label{eq:z2}
Z^{(2)}  =  -\frac14qq_xr^2-q^2rr_x-\frac{i}{4}q^3r^3+\frac{i}{2}qr_{xx},
\end{eqnarray}
\begin{eqnarray}\label{eq:z3}
Z^{(3)}&= & \frac{5}{16}q^4r^4-\frac54q^2r_x^2-\frac32q^2rr_{xx}-\frac32qq_xrr_x\\
& & -\frac14qq_{xx}r^2-2iq^3r^2r_x-\frac34iq^2q_xr^3+\frac{i}{2}qr_{xxx},\nonumber
\end{eqnarray}
and
\begin{eqnarray}\label{eq:z4}
Z^{(4)} & = & 4q^4r^3r_x+\frac{29}{16}q^3q_xr^4\\
& & -\frac92q^2r_xr_{xx}-2q^2rr_{xxx}-\frac{11}{4}qq_xr_x^2-3qq_xrr_{xx}-2qq_{xx}rr_x-\frac14qq_{xxx}r^2\nonumber\\
& & +\frac{7}{16}iq^5r^5\nonumber\\
& & -\frac{15}{4}iq^3r^2r_{xx}-\frac{25}{4}iq^3rr_x^2-8iq^2q_xr^2r_x-iq^2q_{xx}r^3-\frac34iqq_x^2r^3\nonumber\\
& & +\frac{i}{2}qr_{xxxx}.\nonumber
\end{eqnarray}

Let $\phi=Q_1$.
We take
$$ 
q(x)=q_0(x)=q_{1/4,1}(x)=\phi(x)e^{-\frac{i}{2}x+\frac34i\int_{-\infty}^x\phi^2(y)\,dy}.
$$
Then by (\ref{eq:elliptic})
$$
\frac{q_x}{e^{-\frac{i}{2}x+\frac34i\int_{-\infty}^x\phi^2(y)\,dy}}=\phi_x-\frac{i}{2}\phi+\frac34i\phi^3,
$$
$$
\frac{q_{xx}}{e^{-\frac{i}{2}x+\frac34i\int_{-\infty}^x\phi^2(y)\,dy}}=-\frac14\phi+\frac54\phi^3-\frac34\phi^5+i\left(-\phi_x+3\phi^2\phi_x\right),
$$
$$
\frac{q_{xxx}}{e^{-\frac{i}{2}x+\frac34i\int_{-\infty}^x\phi^2(y)\,dy}}=-\frac34\phi_x+6\phi^2\phi_x-6\phi^4\phi_x+i\left(\frac18\phi-\frac{21}{16}\phi^3+3\phi^5-\frac98\phi^7+6\phi\phi_x^2\right),
$$
and
\begin{eqnarray*}
\frac{q_{xxxx}}{e^{-\frac{i}{2}x+\frac34i\int_{-\infty}^x\phi^2(y)\,dy}} &= &\frac{1}{16}\phi-\frac98\phi^3+\frac{45}{8}\phi^5-\frac{111}{16}\phi^7+\frac{63}{32}\phi^9+15\phi\phi_x^2-\frac{57}{2}\phi^3\phi_x^2\\
& &+ i\left(\frac12\phi_x-\frac{15}{2}\phi^2\phi_x+\frac{57}{2}\phi^4\phi_x-\frac{117}{8}\phi^6\phi_x+6\phi_x^3\right).
\end{eqnarray*}
We substitute these into (\ref{eq:z1}), (\ref{eq:z2}), (\ref{eq:z3}) and (\ref{eq:z4}) (or use the recurrence formula (\ref{eq:rec})) to obtain
\begin{eqnarray}\label{eq:z1r}
Z^{(1)}=\frac{1}{2^3}\phi^4-\frac{1}{2^2}\phi^2+\frac{i}{2}\phi\phi_x,
\end{eqnarray}
\begin{eqnarray}\label{eq:z2r}
Z^{(2)}
=\frac{1}{2^2}\phi^3\phi_x-\frac12\phi\phi_x+i\left(-\frac{1}{2^3}\phi^2+\frac{1}{2^2}\phi^4-\frac{1}{2^4}\phi^6\right),
\end{eqnarray}
\begin{eqnarray}\label{eq:z3r}
Z^{(3)} = \frac{1}{2^4}\phi^2-\frac{9}{2^5}\phi^4+\frac{1}{2^3}\phi^6-\frac{1}{2^6}\phi^8+\frac{1}{2^2}\phi^2\phi_x^2+i\left(-\frac{3}{2^3}\phi\phi_x+\frac12\phi^3\phi_x-\frac{1}{2^3}\phi^5\phi_x\right)
\end{eqnarray}
and
\begin{eqnarray}\label{eq:z4r}
Z^{(4)} & = & \frac{1}{2^2}\phi\phi_x-\frac{5}{2^3}\phi^3\phi_x+\frac{5}{2^4}\phi^5\phi_x-\frac{3}{2^6}\phi^7\phi_x+\frac{1}{2^2}\phi\phi_x^3\\
& & +i\left(\frac{1}{2^5}\phi^2-\frac{1}{2^2}\phi^4+\frac{5}{2^5}\phi^6-\frac{5}{2^7}\phi^8+\frac{1}{2^8}\phi^{10}+\frac{5}{2^3}\phi^2\phi_x^2-\frac{3}{2^4}\phi^4\phi_x^2\right).\nonumber
\end{eqnarray}
From (\ref{eq:com}) and the fact that $\phi_x$ is an odd function, it follows that
$$
\int_{\mathbb{R}}Z^{(n)}\,dx=0,\quad n=1,2,3,4,
$$
which imply
$$
P_1(q)=E_1(q)=P_2(q)=E_2(q)=0.
$$

On the other hand, by (\ref{eq:rec}), (\ref{eq:z1r}), (\ref{eq:z2r}), (\ref{eq:z3r}) and (\ref{eq:z4r}), it follows
\begin{eqnarray}\label{eq:z5}
\frac12P_3(q) & = & \Re\int_{\mathbb{R}}Z^{(5)}\,dx\\
& = & \Re\int_{\mathbb{R}}\left(2Z^{(1)}Z^{(3)}+(Z^{(2)})^2+i\phi^2Z^{(4)}+q\frac{\partial}{\partial x}\left(\frac{1}{q}Z^{(4)}\right)\right)\,dx\nonumber\\
& = & \Re\int_{\mathbb{R}}\left(2Z^{(1)}Z^{(3)}+(Z^{(2)})^2+\left(i\phi^2-\frac{q_x}{q}\right)Z^{(4)}\right)\,dx\nonumber\\
& = & \int_{\mathbb{R}}\left(2\Re(Z^{(1)}Z^{(3)})+\Re(Z^{(2)})^2-\frac{\phi_x}{\phi}\Re Z^{(4)}-\left(\frac12+\frac14\phi^2\right)\Im Z^{(4)}\right)\,dx.\nonumber
\end{eqnarray}
From (\ref{eq:z1r}), (\ref{eq:z2r}), (\ref{eq:z3r}) and (\ref{eq:z4r}), we see
\begin{eqnarray*}
\Re(Z^{(1)}Z^{(3)}) =  -\frac{1}{2^6}\phi^4+\frac{5}{2^6}\phi^6-\frac{17}{2^8}\phi^8+\frac{5}{2^8}\phi^{10}-\frac{1}{2^9}\phi^{12}+\frac{3}{2^4}\phi^2\phi_x^2-\frac{5}{2^4}\phi^4\phi_x^2+\frac{3}{2^5}\phi^6\phi_x^2,
\end{eqnarray*}
\begin{eqnarray*}
\Re(Z^{(2)})^2  = -\frac{1}{2^6}\phi^4+\frac{1}{2^4}\phi^6-\frac{5}{2^6}\phi^8+\frac{1}{2^5}\phi^{10}-\frac{1}{2^8}\phi^{12}+\frac{1}{2^2}\phi^2\phi_x^2-\frac{1}{2^2}\phi^4\phi_x^2+\frac{1}{2^4}\phi^6\phi_x^2,
\end{eqnarray*}
\begin{eqnarray*}
\frac{\phi_x}{\phi}\Re Z^{(4)}=\frac14\phi_x^2-\frac{5}{2^3}\phi^2\phi_x^2+\frac{5}{2^4}\phi^4\phi_x^2-\frac{3}{2^6}\phi^6\phi_x^2+\frac{1}{2^2}\phi_x^4,
\end{eqnarray*}
and
\begin{eqnarray*}
\left(\frac12+\frac14\phi^2\right)\Im Z^{(4)} =  \frac{1}{2^6}\phi^2-\frac{15}{2^7}\phi^4+\frac{1}{2^6}\phi^6+\frac{5}{2^8}\phi^8-\frac{1}{2^7}\phi^{10}+\frac{1}{2^{10}}\phi^{12}+\frac{5}{2^4}\phi^2\phi_x^2+\frac{1}{2^4}\phi^4\phi_x^2-\frac{3}{2^6}\phi^6\phi_x^2,
\end{eqnarray*}
which leads to
\begin{eqnarray}\label{eq:zr51}
& & 2\Re(Z^{(1)}Z^{(3)})+\Re(Z^{(2)})^2-\frac{\phi_x}{\phi}\Re Z^{(4)}-\left(\frac12+\frac14\phi^2\right)\Im Z^{(4)}\\
& = & -\frac{1}{2^6}\phi^2+\frac{9}{2^7}\phi^4+\frac{13}{2^6}\phi^6-\frac{59}{2^8}\phi^8+\frac{5}{2^6}\phi^{10}-\frac{9}{2^{10}}\phi^{12}\nonumber\\
& & -\frac{1}{2^2}\phi_x^2+\frac{15}{2^4}\phi^2\phi_x^2-\frac{5}{2^2}\phi^4\phi_x^2+\frac{11}{2^5}\phi^6\phi_x^2-\frac{1}{2^2}\phi_x^4.
\end{eqnarray}
We use (\ref{eq:com}) to have
\begin{eqnarray*}
 \frac{1}{2\pi}P_3(q) 
 & = &  -\frac{1}{2^6}2^2+\frac{9}{2^7}2^3+\frac{13}{2^6}3\cdot 2^3-\frac{59}{2^8}5\cdot 2^4+\frac{5}{2^6}7\cdot 5\cdot 2^3-\frac{9}{2^{10}}9\cdot 7\cdot 2^4\nonumber\\
& & -\frac{1}{2^2}\cdot \frac12+\frac{15}{2^4}-\frac{5}{2^2}\cdot\frac52+\frac{11}{2^5}\cdot 7-\frac{1}{2^2}\frac{3}{2^4}\\
& = & -17+\frac{47}{2}-\frac{13}{2}= 0.
\end{eqnarray*}
Hence $P_3(q)= 0$.

\end{document}